# GRAS–TYPE CONJECTURES
# FOR FUNCTION FIELDS

Cristian D. Popescu

Abstract. Based on results obtained in [15], we construct groups of special $S$–units for function fields of characteristic $p > 0$, and show that they satisfy Gras–type Conjectures. We use these results in order to give a new proof of Chinburg's $\Omega_3$–Conjecture on the Galois module structure of the group of $S$–units, for cyclic extensions of prime degree of function fields.

## 0. Introduction

Let $K/k$ be a finite, abelian extension of function fields of characteristic $p > 0$. Let $G = G(K/k)$ and $g = |G|$. We will denote by $\mathbf{F}_q$ and $\mathbf{F}_{q^\nu}$ the exact fields of constants of $k$ and $K$ respectively, where $q$ is a power of $p$ and $\nu$ is a positive integer. In what follows we will use the same notations as in [15]. For the convenience of the reader, we briefly summarize in this section the main concepts and results of [15] which will be used in our arguments.

For any two finite, nonempty and disjoint sets $S$ and $T$ of primes in $k$, and any field $F$, $k \subseteq F \subseteq K$, $U_{F,S}$ and $U_{F,S,T}$ denote the groups of $S$–units and respectively $(S,T)$–units of $F$; $A_{F,S}$ and $A_{F,S,T}$ are respectively the $S$–ideal class group and $(S,T)$–ideal class group of $F$, as defined in [15, §1.1]. In particular, if $F = K$, we suppress $K$ from the notation, so $U_{K,S} = U_S$, $U_{K,S,T} = U_{S,T}$ etc.

Let us assume for the moment that for a certain positive integer $r$, the set of data $(K/k, S, T, r)$ satisfies the following set of hypotheses:

(H) $\begin{cases} S \neq \emptyset, \quad T \neq \emptyset, \quad S \cap T = \emptyset. \\ S \text{ contains all primes which ramify in } K/k. \\ S \text{ contains at least } r \text{ primes which split completely in } K/k. \\ |S| \geq r + 1. \end{cases}$

Let $(v_1, \ldots, v_r)$ be an ordered $r$–tuple of primes in $S$ which split completely in $K/k$, and let $W = (w_1, \ldots, w_r)$, with $w_i$ prime in $K$, $w_i | v_i$, for every $i = 1, \ldots, r$. One can define a regulator map $\mathbf{C} \overset{r}{\wedge} U_{S,T} \xrightarrow{R_W} \mathbf{C}[G]$, by

$$R_W(u_1 \wedge \cdots \wedge u_r) = \det_{i,j}\left(-\sum_{\sigma \in G} \log |u_j^{\sigma^{-1}}|_{w_i} \cdot \sigma\right), \quad \forall u_1 \wedge \cdots \wedge u_r \in \overset{r}{\wedge} U_{S,T}.$$

Research at MSRI is supported in part by NSF grant DMS–9022140

Typeset by $\mathcal{A}\mathcal{M}\mathcal{S}$-TEX

Throughout this paper, exterior powers are viewed over $\mathbf{Z}[G]$, $|u| = \mathrm{N}w^{-\mathrm{ord}_w(u)}$, for every $u \in K^\times$, and $\mathrm{N}w$ is the order of the residue field $\mathbf{F}_{q^\nu}(w)$ at $w$, for every prime $w$ of $K$. Similarly, if $v$ is a prime in $k$, $\mathrm{N}v$ denotes the order of its residue field $\mathbf{F}_q(v)$.

For every irreducible character $\chi \in \widehat{G}$, let $L_{S,T}(s,\chi)$ be the associated $(S,T)$–L–function, as defined in [15, §1.2], and let $\Theta_{S,T}(s) = \sum_\chi L_{S,T}(s,\chi) e_{\chi^{-1}}$ be the Stickelberger function, where $e_\chi = \sum_{\sigma \in G} \chi(\sigma) \cdot \sigma^{-1} \in \mathbf{C}[G]$. If $r_\chi = \mathrm{ord}_{s=0} L_{S,T}(s,\chi)$, then

$$r_\chi = \begin{cases} \mathrm{card}\{v \in S : \chi|_{G_v} = 1_{G_v}\}, & \text{if } \chi \neq 1_G \\ \mathrm{card}(S) - 1, & \text{if } \chi = 1_G, \end{cases}$$

where $G_v$ is the decomposition group of $v$ in $K/k$. This shows that hypotheses (H) imply that $r_\chi \geq r$, for every $\chi \in \widehat{G}$, and therefore the following definition makes sense in $\mathbf{C}[G]$

$$\Theta_{S,T}^{(r)}(0) \stackrel{\mathrm{def}}{=} \lim_{s \to 0} s^{-r} \Theta_{S,T}(s) = \sum_{\chi \in \widehat{G}} \lim_{s \to 0} s^{-r} L_{S,T}(s,\chi) \cdot e_{\chi^{-1}} \in \mathbf{C}[G].$$

**Definition.** For any ring $R$ with 1, and any $R[G]$–module $M$, let

$$M_{r,S} \stackrel{\mathrm{def}}{=} \{m \in M \,|\, e_\chi \cdot m = 0 \text{ in } M \otimes_R \mathbf{C}, \, \forall \chi \in \widehat{G} \text{ with } r_\chi \neq r\}.$$

In [15] (see Theorem 3.2.1 and Corollary 3.2.2) we proved the following result:

**Theorem 0.** *If the set of data $(K/k, S, T, r)$ satisfies hypotheses (H) then:*

(1) *There exists a unique $\varepsilon_{S,T} \in \left(\mathbf{Z}[1/g] \stackrel{r}{\wedge} U_{S,T}\right)_{r,S}$ such that*

$$R_W(\varepsilon_{S,T}) = \Theta_{S,T}^{(r)}(0).$$

(2) *The element $\varepsilon_{S,T}$ satisfies the equality:*

$$\mathbf{Z}[1/g][G]\,\varepsilon_{S,T} = \mathrm{Fitt}_{\mathbf{Z}[G]}(A_{S,T}) \cdot \left(\mathbf{Z}[1/g] \stackrel{r}{\wedge} U_{S,T}\right)_{r,S}.$$

Throughout this paper, for any noetherian ring $R$ with 1 and any finitely generated $R$–module $M$, $\mathrm{Fitt}_R(M)$ denotes the Fitting ideal of $M$. For the definition and properties of the Fitting ideals needed for our purposes, the reader can consult [15, §1.4].

*Let us now fix once and for all two finite, nonempty, disjoint sets $S$ and $T$ of primes in $k$, $S$ containing at least one prime which splits completely in $K/k$, and $T$ not containing primes which ramify in $K/k$.*

We will use Theorem 0(1) in order to construct two $\mathbf{Z}[G]$–modules of special units $\mathcal{E}_{S,T}$ and $\mathcal{E}_S$, of finite index in $U_{S,T}$ and $U_S$ respectively. Theorem 0(2) will help us prove that statements similar to the ones conjectured by Gras in [9] and proved by Mazur and Wiles in [12] for the classical number field case of cyclotomic units, are satisfied by $\mathcal{E}_{S,T}$ and $\mathcal{E}_S$ (Theorems 1.4, 2.2, 3.10). We use these results in §4 in order to give a new proof to a special case of Chinburg's $\Omega_3$–Conjecture (Theorem 4.2.7).

Before proceeding we would like to make a final useful remark on group–rings and their modules (see also [15, §1.3]). If $L$ is a field of characteristic 0, then $\widehat{G}(L)$ will denote the set of characters associated to $L$–irreducible representations of $G$. If $\overline{L}$ is an algebraic closure of $L$, then $G(\overline{L}/L)$ acts canonically on $\widehat{G}(\overline{L})$ and $\widehat{G}(L)$ can be viewed as the set of orbits with respect to this action. For $\chi \in \widehat{G}(\overline{L})$ and $\psi \in \widehat{G}(L)$, we write $\chi | \psi$ if $\chi$ is in the orbit represented by $\psi$. For a subgroup $H \subseteq G$ (or an element $\sigma \in G$) we write $\psi|_H = \mathbf{1}_H$ (or $\psi(\sigma) = 1$) if $\chi|_H = \mathbf{1}_H$ (or $\chi(\sigma) = 1$), for some (i.e. all) $\chi \in \widehat{G}(\overline{L})$, such that $\chi|\psi$.

If $R$ is a Dedekind domain containg $\mathbf{Z}[1/g]$, and $L$ is its field of fractions, then one has a canonical decomposition

$$R[G] = \bigoplus_{\psi \in \widehat{G}(L)} D_\psi,$$

where $D_\psi = R[G] \cdot e_\psi$ are finite extensions of $R$, and $e_\psi = 1/g \sum_{\sigma \in G} \psi(\sigma) \cdot \sigma^{-1} \in R[G]$. If $M$ is an $R[G]$–module, then one has a decomposition

$$M = \bigoplus_{\psi \in \widehat{G}(L)} M^\psi,$$

where $M^\psi = M \underset{R[G]}{\otimes} D_\psi$, for every $\psi \in \widehat{G}(L)$.

**Acknowledgements.** It is a pleasure to thank Karl Rubin for his invaluable help and support, and Adebisi Agboola for helpful conversations leading to the result in §4.

## 1. The Group $\mathcal{E}_{S,T}$ of Stark $(S,T)$–units

Throughout this section the sets of primes $S$ and $T$ as above will be fixed, so that $S$ contains at least one prime which splits completely and $T$ does not contain primes which ramify in $K/k$. Let $F$ be an intermediate field of $K/k$ and let $S'$ be a finite set of primes in $k$ containing $S$. We define

$$r_{F,S'} = \begin{cases} \text{card}\{v \in S' \,|\, v \text{ splits completely in } F/k\}, & \text{for } F \neq k \\ \text{card}\,(S') - 1, & \text{for } F = k. \end{cases}$$

In particular, if $S' = S$, we will make the notation $r_{F,S} = r_F$.

Let $L$ be any finite field of characteristic 0, and let $\psi \in \widehat{G}(L)$. We define

$$r_{\psi,S'} = \begin{cases} \text{card}\{v \in S' \,|\, \psi|_{G_v} = \mathbf{1}_{G_v}\}, & \text{for } \psi \neq \mathbf{1}_G \\ \text{card}\,(S') - 1, & \text{for } \psi = \mathbf{1}_G. \end{cases}$$

In particular, for $S = S'$, we make the notation $r_{\psi,S} = r_\psi$.

For $\psi \in \widehat{G}(L)$ as above, let $H_\psi = \ker(\psi)$ and let $K_\psi$ be the maximal subfield of $K$ fixed by $H_\psi$. The injectivity of $\psi$ on $G(K_\psi/k)$ obviously implies that $r_{K_\psi,S'} = r_{\psi,S'}$, for any $S'$ as above.

Let $S = \{v_0, \ldots, v_s\}$. For every $i = 0, \ldots, s$, let us fix a prime $w_i$ above $v_i$ in $K$ and let $w_{i,F}$ denote the prime in $F$ sitting below $w_i$, for every intermediate field $F$. Let $W_F$ be the ordered $r_F$–tuple of primes in $F$ defined by

$$W_F = \begin{cases} (w_{i,F} \,|\, v_i \text{ splits completely in } F), & \text{for } F \neq k \\ (v_0, \ldots, v_{s-1}), & \text{for } F = k, \end{cases}$$

ordered so that $w_{i,F}$ preceeds $w_{j,F}$ iff $i < j$.

**Definition 1.1.** A pair $(F, S')$ consisting of an intermediate field $F$, $F \neq k$, and a finite set of primes $S'$ in $k$ is called $(S,T)$–admissible if the following conditions are satisfied:

(1) $S \subseteq S'$
(2) $|S'| \geq |S| + 2$
(3) the set of data $(F/k, S', T, r_F)$ satisfies hypotheses (H)

The pair $(k, S)$ is also declared to be $(S,T)$–admissible.

Let us emphasize that if $(F, S')$ is an $(S,T)$–admissible pair, then the set of data $(F/k, S', T, r_F)$ satisfies hypotheses (H). For any such pair $(F, S')$, let

$$\varepsilon_{F,S',T} \in \left[\mathbf{Z}\,[1/g] \overset{r_F}{\wedge} U_{F,S',T}\right]_{r_{F,S'}}$$

be the unique element associated by Theorem 0 to the set of data $(F/k, S', T, r_F)$ and the regulator $R_{W_F}$.

Let $U^*_{F,S',T} = \text{Hom}_{\mathbf{Z}[G]}(U_{F,S',T}, \mathbf{Z}\,[G])$ be the $\mathbf{Z}\,[G]$–dual of $U_{F,S',T}$. Then for every $\Phi = \phi_1 \wedge \cdots \wedge \phi_{r_F - 1} \in \mathbf{Z}\,[1/g] \overset{r_F - 1}{\wedge} U^*_{F,S',T}$ we get a $\mathbf{Z}\,[1/g]\,[G]$–morphism

$$\mathbf{Z}\,[1/g] \overset{r_F}{\wedge} U_{F,S',T} \xrightarrow{\Phi} \mathbf{Z}\,[1/g]\, U_{F,S',T},$$

given by
$$\Phi(u_1 \wedge \cdots \wedge u_{r_F}) = \sum_{1 \leq i \leq r_F} (-1)^i \det_{\substack{1 \leq k,j \leq r_F \\ j \neq i}} (\phi_k(u_j)) \cdot u_i,$$

for all $u_1 \wedge \cdots \wedge u_{r_F} \in \bigwedge^{r_F} U_{F,S',T}$.

For every intermediate field $F$, let $\mu_F$ be the group of roots of unity in $F$, and let $e_F = |\mu_F|$. In particular, if $F = K$, we have $\mu_K = \mathbf{F}_{q^\nu}^\times$ and therefore $e_K = q^\nu - 1$.

**Definition 1.2.** Let $\widetilde{\mathcal{E}}_{S,T}$ be the $\mathbf{Z}[1/g][G]$–submodule of $\mathbf{Z}[1/g] \bigotimes_{\mathbf{Z}} (K^\times/\mu_K)$ generated by the set

$$\{\Phi(\varepsilon_{F,S',T}) \mid (F,S') \text{ is } (S,T)\text{–admissible}, \Phi \in \mathbf{Z}[1/g]^{r_F-1} \wedge U^*_{F,S',T}\}.$$

**Definition 1.3.** The group $\mathcal{E}_{S,T}$ of Stark $(S,T)$–units is defined by

$$\mathcal{E}_{S,T} = \widetilde{\mathcal{E}}_{S,T} \bigcap U_{S,T},$$

where the intersection is viewed inside $\mathbf{Z}[1/g] \bigotimes_{\mathbf{Z}} (K^\times/\mu_K)$.

The main goal of this section is the proof of the following Gras–type Conjecture for the group $\mathcal{E}_{S,T}$ of special $(S,T)$–units:

**Theorem 1.4.** *Let $K/k$, $S$ and $T$ be as above. Then*
  (1) *The index $[U_{S,T} : \mathcal{E}_{S,T}]$ is finite.*
  (2) *For any prime number $\ell$, such that $\gcd(\ell, g) = 1$, and any $\psi \in \widehat{G}(\mathbf{Q}_\ell)$*

$$|(U_{S,T}/\mathcal{E}_{S,T} \otimes \mathbf{Z}_\ell)^\psi| = |(A_{S,T} \otimes \mathbf{Z}_\ell)^\psi|^{r_\psi}.$$

Before proceeding to the proof of the theorem above, we need a few more lemmas and definitions.

**Definition 1.5.** Let $L$ be a field of characteristic 0, and let $\psi \in \widehat{G}(L)$. A pair $(F, S')$ is called $(S, T, \psi)$–admissible if the following conditions are satisfied
  (1) $(F, S')$ is $(S, T)$–admissible
  (2) $K_\psi \subseteq F$
  (3) $r_{\psi, S'} = r_F$.

Let us notice that $(k, S)$ is the only $(S, T, \mathbf{1}_G)$–admissible pair, while for $\psi \neq \mathbf{1}_G$ there exist infinitely many $(S, T, \psi)$–admissible pairs.

**Lemma 1.6.**
(1) Let $\psi \in \widehat{G}(L)$ and let $(F, S')$ be an $(S, T, \psi)$–admissible pair. Then $r_{K_\psi} = r_{\psi, S'} = r_F$ and $(K_\psi, S')$ is also $(S, T, \psi)$–admissible.
(2) If $(F, S')$ and $(F, S'')$ are $(S, T, \psi)$–admissible pairs, then $(F, S' \bigcup S'')$ is $(S, T, \psi)$–admissible as well.

**Proof.** Both statements are clear if $\psi = \mathbf{1}_G$. Let $\psi \neq \mathbf{1}_G$.
(1) From the definitions above we have the following sequence of inequalities

$$r_{K_\psi} \geq r_F = r_{\psi, S'} \geq r_\psi = r_{K_\psi},$$

which shows that $r_{K_\psi} = r_{\psi, S'} = r_F$. This obviously implies that $(K_\psi, S')$ is $(S, T, \psi)$–admissible.
(2) Since $r_F = r_{\psi, S'} = r_{\psi, S''}$, we have $\psi|_{G_v} \neq \mathbf{1}_{G_v}$, for every $v \in S' \setminus S$ and every $v \in S'' \setminus S$. This shows that $r_F = r_{\psi, S' \cup S''}$, which obviously implies that $(F, S' \bigcup S'')$ is $(S, T, \psi)$–admissible. $\square$

**Lemma 1.7.** Let $\psi \in \widehat{G}(\mathbf{Q})$ and let $(F, S')$ be an $(S, T, \psi)$–admissible pair. Then the natural $\mathbf{Z}[G]$–morphisms involved in the following commutative diagrams:

$$\begin{array}{ccc} U_{S,T} \longrightarrow U_{S',T} & & A_{S,T} \longrightarrow A_{S',T} \\ \uparrow \quad\quad \uparrow & & \uparrow \quad\quad \uparrow \\ U_{F,S,T} \longrightarrow U_{F,S',T} & & A_{F,S,T} \longrightarrow A_{F,S',T}, \end{array}$$

become $\mathbf{Z}[1/g][G]^\psi$–isomorphisms, when tensored with $\mathbf{Z}[1/g][G]^\psi$ over $\mathbf{Z}[G]$.

**Proof.** The vertical morphisms in both diagrams obviously become isomorphisms when tensored with $\mathbf{Z}[1/g][G]^\psi$. (Their inverses are the $\mathbf{Z}[1/g][G]^\psi$–linear extensions of the norm map from $K$ to $F$ at the level of units and ideal class–groups respectively.) This observation settles the lemma above for $\psi = \mathbf{1}_G$ and it shows that it is enough to prove that the lower horizontal morphisms in both diagrams become isomorphisms when tensored with $\mathbf{Z}[1/g][G]^\psi$, for $\psi \neq \mathbf{1}_G$.

Let $\psi \neq \mathbf{1}_G$ and let $\mathcal{S}$ be the $\mathbf{Z}[G]$–submodule of $A_{F,S,T}$ generated by the ideal classes

$$\{w \mid w \text{ prime in } F, w|v \text{ for some } v \in S' \setminus S\}.$$

Let us fix a prime $w(v)$ in $F$ for every prime $v \in S' \setminus S$. We have two exact sequences of $\mathbf{Z}[G]$–modules (see [15, §2 (9)] and [17, §5.1]):

$$0 \longrightarrow \mathcal{S} \longrightarrow A_{F,S,T} \longrightarrow A_{F,S',T} \longrightarrow 0$$

$$0 \longrightarrow U_{F,S,T} \longrightarrow U_{F,S',T} \longrightarrow \bigoplus_{v \in S' \setminus S} \mathbf{Z}[G] w(v) \longrightarrow \mathcal{S} \longrightarrow 0.$$

Let us fix $v \in S' \setminus S$. Since $\psi \neq \mathbf{1}_G$ and $(F, S')$ is $(S, T, \psi)$–admissible, we have $\psi|_{G_v} \neq \mathbf{1}_{G_v}$. This shows that there exists $\sigma \in G$ such that $^\sigma w(v) = w(v)$ and $\psi(\sigma) \neq 1$. We therefore have an equality

$$e_\psi(\sigma - 1) \cdot w(v) = 0 \tag{1}$$

in $\mathbf{Z}[1/g][G]^\psi$, where $e_\psi$ is the idempotent associated to $\psi$ in $\mathbf{Z}[1/g][G]$. The fact that $\psi(\sigma) \neq 1$ easily implies that $e_\psi(\sigma - 1)$ is invertible in $\mathbf{Z}[1/g][G]^\psi$ and therefore relation (1) above shows that

$$\bigoplus_{v \in S' \setminus S} \mathbf{Z}[1/g][G]^\psi \cdot w(v) = 0.$$

The last equality proves that, when tensored with $\mathbf{Z}[1/g][G]^\psi$ (which is a flat $\mathbf{Z}[G]$–algebra), the exact sequences above yield isomorphisms

$$(\mathbf{Z}[1/g] A_{F,S,T})^\psi \xrightarrow{\sim} (\mathbf{Z}[1/g] A_{F,S',T})^\psi, \quad (\mathbf{Z}[1/g] U_{F,S,T})^\psi \xrightarrow{\sim} (\mathbf{Z}[1/g] U_{F,S',T})^\psi,$$

which concludes the proof of Lemma 1.7. $\square$

**Lemma 1.8.**
(1) Let $\psi \in \widehat{G}(\mathbf{Q})$. Then $\widetilde{\mathcal{E}}_{S,T}^\psi$ is generated as a $\mathbf{Z}[1/g][G]^\psi$–module by the set

$$\{e_\psi \cdot \Phi\left(\varepsilon_{K_\psi, S', T}\right) \mid \Phi \in \mathbf{Z}[1/g] \stackrel{r_\psi - 1}{\wedge} U^*_{K_\psi, S', T}\},$$

for any $S'$ such that $(K_\psi, S')$ is $(S, T, \psi)$–admissible.

(2) $\widetilde{\mathcal{E}}_{S,T}^\psi \subseteq \mathbf{Z}[1/g] U_{S,T}$ and $\mathbf{Z}[1/g] \mathcal{E}_{S,T} = \widetilde{\mathcal{E}}_{S,T}$.

**Proof.** (1) Let $(F, S')$ be an $(S, T)$–admissible pair. If the pair $(F, S')$ is not $(S, T, \psi)$–admissible, then either $K_\psi \not\subseteq F$, or $K_\psi \subseteq F$ and $r_{\psi, S'} \neq r_F$. In both circumstances $\varepsilon_{F,S',T}^\psi = 0$ in $\left(\mathbf{Z}[1/g] \stackrel{r_F}{\wedge} U_{F,S',T}\right)^\psi$. We thus have

$$\Phi\left(\varepsilon_{F,S',T}^\psi\right) = 0, \text{ for all } \Phi \in \mathbf{Z}[1/g] \stackrel{r_F - 1}{\wedge} U^*_{F,S',T}.$$

It is therefore enough to restrict ourselves to $(S, T, \psi)$–admissible pairs in the definition of $\widetilde{\mathcal{E}}_{S,T}^\psi$. Let $(F, S')$ be such a pair. Then, according to Lemma 1.6, $(K_\psi, S')$ is also $(S, T, \psi)$–admissible and $r_F = r_{K_\psi} = r_\psi$. The definition and uniqueness of $\varepsilon_{K_\psi, S', T}$ and $\varepsilon_{F, S', T}$ show that

$$\varepsilon_{F,S',T} = \mathrm{N}^{(r_\psi)}_{F/K_\psi}\left(\varepsilon_{K_\psi, S', T}\right), \tag{2}$$

where $\mathrm{N}_{F/K_\psi}^{(r_\psi)}$ is the $\mathbf{Z}[1/g]$–linear extension to $\mathbf{Z}[1/g] \overset{r_\psi}{\wedge} U_{F,S',T}$ of the $r_\psi$–exterior power of the norm map

$$U_{F,S',T} \xrightarrow{\mathrm{N}_{F/K_\psi}} U_{K_\psi,S',T} \ .$$

(see [17, §6] for (2)). The $\mathbf{Z}$–freeness of $U_{F,S',T}$ gives an isomorphism of abelian groups

$$U^*_{F,S',T} \xrightarrow[\sim]{\phi \to \phi_0} \mathrm{Hom}_{\mathbf{Z}}(U_{F,S',T}, \mathbf{Z}) \ , \tag{3}$$

defined by $\phi(u) = \sum_{\sigma \in G(F/k)} \phi_0(\sigma^{-1}u) \cdot \sigma$, for every $u \in U_{F,S',T}$. Obviously the same type of isomorphism holds true for $U^*_{K_\psi,S',T}$ as well. The inclusion $U_{K_\psi,S',T} \subseteq U_{F,S',T}$ induces therefore a surjective morphism

$$\mathbf{Z}[1/g] U^*_{F,S',T} \xrightarrow{\phi \to \tilde\phi} \mathbf{Z}[1/g] U^*_{K_\psi,S',T} ,$$

satisfying the relation $\tilde\phi(\mathrm{N}_{F/K_\psi} u) = \pi \phi(u)$, for every $u \in U_{K_\psi,S',T}$, where $\pi$ is the natural projection

$$\mathbf{Z}[G(F/k)] \xrightarrow{\pi} \mathbf{Z}[G(K_\psi/k)] \ .$$

We therefore obtain a surjective $\mathbf{Z}[1/g]$–morphism

$$\mathbf{Z}[1/g] \overset{r_\psi - 1}{\wedge} U^*_{F,S',T} \longrightarrow \mathbf{Z}[1/g] \overset{r_\psi - 1}{\wedge} U^*_{K_\psi,S',T}$$

$$\Phi = \phi_1 \wedge \cdots \wedge \phi_{r_\psi - 1} \longrightarrow \tilde\Phi = \tilde\phi_1 \wedge \cdots \wedge \tilde\phi_{r_\psi - 1}$$

satisfying the relation

$$\tilde\Phi\left(\mathrm{N}_{F/K_\psi}^{(r_\psi)}(u)\right) = \mathrm{N}_{F/K_\psi}(\Phi(u)) \ ,$$

for every $u \in \mathbf{Z}[1/g] \overset{r_\psi}{\wedge} U_{F,S',T}$. The equality above obviously implies that, for any $(S,T,\psi)$–admissible pair $(F,S')$ and any $\Phi \in \mathbf{Z}[1/g] \overset{r_\psi - 1}{\wedge} U^*_{F,S',T}$, the following relation holds true in $\widetilde{\mathcal{E}}^\psi_{S,T}$:

$$[F : K_\psi] e_\psi \Phi(\varepsilon_{F,S',T}) = e_\psi \tilde\Phi(\varepsilon_{K_\psi,S',T}) \ . \tag{4}$$

Let now $S'$ and $S''$ be two sets of primes such that the pairs $(K_\psi, S')$ and $(K_\psi, S'')$ are both $(S,T,\psi)$–admissible. According to Lemma 1.6(2), the pair $(K_\psi, S' \cup S'')$ is $(S,T,\psi)$–admissible as well. The definition and uniqueness of the elements $\varepsilon$ give the following relations (see [17, §6]):

$$\varepsilon_{K_\psi, S' \cup S'', T} = \left[\prod_{v \in S'' \setminus S'} (1 - \sigma_v^{-1})\right] \varepsilon_{K_\psi, S', T} = \left[\prod_{v \in S' \setminus S''} (1 - \sigma_v^{-1})\right] \varepsilon_{K_\psi, S'', T} ,$$

where $\sigma_v$ is the Frobenius morphism associated to $v$ in $G(K_\psi/k)$. This implies that

$$\left[\prod_{v \in S'' \setminus S'} e_\psi \left(1 - \sigma_v^{-1}\right)\right] \varepsilon_{K_\psi, S', T} = \left[\prod_{v \in S' \setminus S''} e_\psi \left(1 - \sigma_v^{-1}\right)\right] \varepsilon_{K_\psi, S'', T}. \quad (5)$$

Relations (4) and (5), the fact that $[F : K_\psi]$ and $e_\psi(\sigma_v - 1)$ are invertible in $\mathbf{Z}[1/g][G]^\psi$, for $v \in (S' \setminus S'') \cup (S'' \setminus S')$ (see the proof of Lemma 1.7), and Lemma 1.7 show that (1) holds true.

(2) Lemma 1.7 implies that for any $S'$ such that $(K_\psi, S')$ is $(S, T, \psi)$–admissible, we have $\mathbf{Z}[1/g][G]^\psi U_{K_\psi, S', T} = \mathbf{Z}[1/g][G]^\psi U_{S, T}$. We thus have

$$e_\psi \Phi\left(\varepsilon_{K_\psi, S', T}\right) \in \mathbf{Z}[1/g][G]^\psi U_{S, T},$$

for any $\Phi \in \overset{r_\psi - 1}{\wedge} U^*_{K_\psi, S', T}$. According to this relation, the fact that $\widetilde{\mathcal{E}_{S,T}} \subseteq \mathbf{Z}[1/g] U_{S, T}$ follows from statement (1) in our lemma. The equality $\mathbf{Z}[1/g]\mathcal{E}_{S,T} = \widetilde{\mathcal{E}_{S,T}}$ is now an obvious consequence of the definition of $\mathcal{E}_{S,T}$. $\square$

**Proof of Theorem 1.4.** Since $U_{S,T}$ is a finitely generated $\mathbf{Z}$–module, it is obviously enough to prove the second statement of the theorem. Let $\ell$ be a prime number such that $\gcd(\ell, g) = 1$, and let

$$\mathbf{Z}_\ell[G] = \bigoplus_{\psi \in \widehat{G}(\mathbf{Q}_\ell)} D_\psi$$

be the decomposition of $\mathbf{Z}_\ell[G]$ into a direct sum of discrete valuation rings, as described in §0. Let us fix $\psi \in \widehat{G}(\mathbf{Q}_\ell)$ and let $S'$ be a set of primes such that $(K_\psi, S')$ is $(S, T, \psi)$–admissible. Since $r_{\psi, S'} = r_\psi$, Lemma 1.5.1(1) of [15] gives us the following isomorphisms of $D_\psi$–modules

$$\left(\mathbf{Z}_\ell U_{K_\psi, S', T}\right)^\psi \xrightarrow{\sim} D_\psi^{r_\psi}, \quad \mathrm{Hom}_{D_\psi}\left(\left(\mathbf{Z}_\ell U_{K_\psi, S', T}\right)^\psi, D_\psi\right) \xrightarrow{\sim} D_\psi^{r_\psi}.$$

Let $\{e_1, \ldots, e_{r_\psi}\}$ be a $D_\psi$ basis of $\left(\mathbf{Z}_\ell U_{K_\psi, S', T}\right)^\psi$ and let $\{e_1^*, \ldots, e_{r_\psi}^*\}$ be its $D_\psi$–dual. We have the following equalities of $D_\psi$–modules

$$D_\psi \bigotimes_{\mathbf{Z}[1/g][G]} \left(\mathbf{Z}[1/g][G] \overset{r_\psi}{\wedge} U_{K_\psi, S', T}\right)_{r_\psi, S'} = \left(\mathbf{Z}_\ell \overset{r_\psi}{\wedge} U_{K_\psi, S', T}\right)^\psi$$
$$= D_\psi \left(e_1 \wedge \cdots \wedge e_{r_\psi}\right).$$

Let us now combine Theorem 0(2), applied to the set of data $(K_\psi/k, S', T, r_\psi)$, with the equalities above. After tensoring with $D_\psi$ over $\mathbf{Z}[1/g][G]$, we obtain

$$D_\psi \varepsilon_{K_\psi, S', T} = a_\psi \left(e_1 \wedge \cdots \wedge e_{r_\psi}\right), \quad (6)$$

where

$$a_\psi = \mathrm{Fitt}_{\mathbf{Z}[G]}\left(A_{K_\psi,S',T}\right)\bigotimes_{\mathbf{Z}[G]} D_\psi = \mathrm{Fitt}_{D_\psi}\left(A_{K_\psi,S',T}\bigotimes_{\mathbf{Z}[G]} D_\psi\right). \quad (7)$$

(For the last equality above see [15, §1.4(b)].) We are now prepared to compute the $D_\psi$–module $(\mathbf{Z}_\ell \mathcal{E}_{S,T})^\psi$ in terms of the chosen basis $\{e_1, \ldots, e_{r_\psi}\}$. Let us first notice that (3) shows that we have the following relations:

$$\left(\mathbf{Z}_\ell U^*_{K_\psi,S',T}\right) \xrightarrow{\sim} \mathrm{Hom}_{\mathbf{Z}_\ell[G]}\left(\mathbf{Z}_\ell U_{K_\psi,S',T}, \mathbf{Z}_\ell[G]\right) =$$
$$= \bigoplus_{\chi \in \widehat{G}(\mathbf{Q}_\ell)} \mathrm{Hom}_{D_\chi}\left((\mathbf{Z}_\ell U_{K_\psi,S',T})^\chi, D_\chi\right).$$

This direct sum decomposition together with the definition of $\widetilde{\mathcal{E}}_{S,T}$, Lemma 1.8 and equality (6) imply that

$$\begin{aligned}
(\mathbf{Z}_\ell \mathcal{E}_{S,T})^\psi = \left(\mathbf{Z}_\ell \widetilde{\mathcal{E}}_{S,T}\right)^\psi = \\
= \sum_{1 \le i \le r_\psi} \left(e_1^* \wedge \cdots \wedge \hat{e}_i^* \wedge \cdots \wedge e_{r_\psi}^*\right)\left(D_\psi \varepsilon_{K_\psi,S',T}\right) = \\
= \sum_{1 \le i \le r_\psi} \left(e_1^* \wedge \cdots \wedge \hat{e}_i^* \wedge \cdots \wedge e_{r_\psi}^*\right)\left(a_\psi \cdot e_1 \wedge \cdots \wedge e_{r_\psi}\right) = \\
= a_\psi e_1 \oplus \cdots \oplus a_\psi e_{r_\psi},
\end{aligned} \quad (8)$$

equality viewed inside $\left(\mathbf{Z}_\ell U_{K_\psi,S',T}\right)^\psi = (\mathbf{Z}_\ell U_{S,T})^\psi = D_\psi^{r_\psi}$ (see Lemma 1.7 above and [15], Lemma 1.5.1(1)). Relation (8), together with (7) and [15, §1.4(f)], obviously give us the following equality

$$\left|\left(U_{S,T}/\mathcal{E}_{S,T} \otimes \mathbf{Z}_\ell\right)^\psi\right| = [D_\psi : a_\psi]^{r_\psi} = \left|\left(A_{K_\psi,S',T} \otimes \mathbf{Z}_\ell\right)^\psi\right|^{r_\psi}.$$

On the other hand, after tensoring with $D_\psi$, Lemma 1.7 shows that

$$\left(A_{K_\psi,S',T} \otimes \mathbf{Z}_\ell\right)^\psi \xrightarrow{\sim} (A_{S,T} \otimes \mathbf{Z}_\ell)^\psi,$$

as $D_\psi$–modules, which, combined with the previous equality, settles statement (2) in our theorem. $\square$

## 2. The Group $\mathcal{E}_S$ of Stark $S$–units

With the same assumptions as in the previous section, let us keep $S$ fixed and let $T$ vary, so that $S \cap T = \emptyset$, and $T$ does not contain any prime which ramifies in $K/k$. We will denote by $\mathcal{T}$ the set of all possible choices of sets $T$. For each $T \in \mathcal{T}$ we have just defined a $\mathbf{Z}[G]$–submodule $\mathcal{E}_{S,T}$ of $U_{S,T}$.

**Definition 2.1.** The group of Stark $S$–units $\mathcal{E}_S$ is the $\mathbf{Z}[G]$–submodule of $U_S$, generated by

$$\mu_K \bigcup \left( \bigcup_{T \in \mathcal{T}} \mathcal{E}_{S,T} \right).$$

The main goal of this section is the proof of the following Gras–type Conjecture satisfied by the group $\mathcal{E}_S$ of special $S$–units:

**Theorem 2.2.** *Let $S$ be a finite set of primes in $k$, containing at least one prime which splits completely in $K/k$. Then:*
  (1) *The index $[U_S : \mathcal{E}_S]$ is finite.*
  (2) *For any prime number $\ell$, such that $\gcd(\ell, g \cdot e_K) = 1$, and any character $\psi \in \widehat{G}(\mathbf{Q}_\ell)$,*

$$|(U_S/\mathcal{E}_S \otimes \mathbf{Z}_\ell)^\psi| = |(A_S \otimes \mathbf{Z}_\ell)^\psi|^{r_\psi}.$$

In order to prove the theorem above, we need the following

**Lemma 2.3.** *For any $T \in \mathcal{T}$, $|\mathcal{E}_S/\mathcal{E}_{S,T}|$ is divisible only by primes dividing $g \cdot \prod_{w \in T_K} (Nw - 1)$.*

**Proof.** Let $T_1, T_2 \in \mathcal{T}$ and let $\psi \in \widehat{G}(\mathbf{Q})$. Let $S'$ be a finite set of primes in $k$ such that the pair $(K_\psi, S')$ is $(S, T_1 \cup T_2, \psi)$–admissible. Then $(K_\psi, S')$ is obviously $(S, T_1, \psi)$–admissible and $(S, T_2, \psi)$–admissible as well. The definition and uniqueness of the elements $\varepsilon$ associated by Theorem 0 to the sets of data $(K_\psi, S', T_i, r_\psi)$, for $i = 1, 2$, and $(K_\psi, S', T_1 \cup T_2, r_\psi)$, give the following relations:

$$\varepsilon_{K_\psi, S', T_1 \cup T_2} = \left[ \prod_{v \in T_2 \setminus T_1} (1 - \sigma_v^{-1} \cdot Nv) \right] \cdot \varepsilon_{K_\psi, S', T_1} = \\ = \left[ \prod_{v \in T_1 \setminus T_2} (1 - \sigma_v^{-1} \cdot Nv) \right] \cdot \varepsilon_{K_\psi, S', T_2} \qquad (9)$$

(see [17, §6] for (9)). Let us also notice (see Lemma 3.1 below as well) that for any prime $v$ in $k$, which is unramified in $K/k$, we have

$$\left(1 - \sigma_v^{-1} \cdot Nv\right) \mu_K = \{1\}, \qquad (10)$$

where $\sigma_v$ is the Frobenius morphism associated to $v$ in $G(K/k)$. In light of Lemma 1.8(1), equalities (9) and (10) imply that, for any set $T \in \mathcal{T}$,

$$\prod_{v \in T} (1 - \sigma_v^{-1} \cdot Nv) \cdot \mathbf{Z}[1/g]\mathcal{E}_S \subseteq \mathbf{Z}[1/g]\mathcal{E}_{S,T},$$

and therefore

$$\prod_{w \in T_K} (1 - Nw) \cdot \mathbf{Z}[1/g]\mathcal{E}_S \subseteq \mathbf{Z}[1/g]\mathcal{E}_{S,T}. \tag{11}$$

(We are using the fact that $(1 - \sigma_v^{-1} \cdot Nv)$ divides $(1 - Nw)$ in the group ring $\mathbf{Z}[G]$, for any $v \in T$, and any $w$ in $K$, such that $w|v$.) The statement in the lemma follows directly from (11). $\square$

**Proof of Theorem 2.2.** Since $U_S$ is a finitely generated $\mathbf{Z}$–module, it is obviously enough to prove (2) in the statement of the Theorem above.

Due to the fact that $K$ has a divisor of degree 1 over $\mathbf{F}_{q^\nu}$, supported outside any finite set of primes (see [13]), we have

$$\gcd\{(1 - Nw); w \in T_K, T \in \mathcal{T}\} = e_K.$$

Let us fix a prime number $\ell$, such that $\gcd(\ell, ge_K) = 1$, and let us consider $T \in \mathcal{T}$, $T = \{v\}$, such that $\gcd(1 - Nw, \ell) = 1$, for every prime $w \in T_K$. We have the following exact sequences of $\mathbf{Z}[G]$–modules (see [15, §1.1, (1)]):

$$0 \longrightarrow U_S/U_{S,T} \longrightarrow \bigoplus_{w \in T_K} \mathbf{F}_{q^\nu}(w)^\times \longrightarrow A_{S,T} \longrightarrow A_S \longrightarrow 0$$

$$0 \longrightarrow (\mathcal{E}_S \cap U_{S,T})/\mathcal{E}_{S,T} \longrightarrow U_{S,T}/\mathcal{E}_{S,T} \longrightarrow U_S/\mathcal{E}_S \longrightarrow U_S/(U_{S,T} \cdot \mathcal{E}_S) \longrightarrow 0.$$

Lemma 2.3 implies that, when tensored with $\mathbf{Z}_\ell$, the exact sequences above give the following isomorphisms of $\mathbf{Z}_\ell[G]$–modules

$$A_{S,T} \otimes \mathbf{Z}_\ell \xrightarrow{\sim} A_S \otimes \mathbf{Z}_\ell, \qquad U_{S,T}/\mathcal{E}_{S,T} \otimes \mathbf{Z}_\ell \xrightarrow{\sim} U_S/\mathcal{E}_S \otimes \mathbf{Z}_\ell.$$

These isomorphisms, combined with Theorem 1.4(2), imply that

$$|(U_S/\mathcal{E}_S \otimes \mathbf{Z}_\ell)^\psi|^{r_\psi} = |(U_{S,T}/\mathcal{E}_{S,T} \otimes \mathbf{Z}_\ell)^\psi|^{r_\psi} =$$
$$= |(A_{S,T} \otimes \mathbf{Z}_\ell)^\psi|^{r_\psi} = |(A_S \otimes \mathbf{Z}_\ell)^\psi|^{r_\psi},$$

for every character $\psi \in \widehat{G}(\mathbf{Q}_\ell)$. This concludes the proof of Theorem 2.2. $\square$

**3. The case $|S| = 1$**

We keep the same notations as in the previous sections, and we assume in addition that $S$ consists of a single prime $v_0$ which splits completely in $K/k$. (If one thinks of $v_0$ as the prime at infinity, then this case is the function field analogue of the totally real number field case.) The goal of this section is to remove condition $\gcd(\ell, e_K) = 1$ in Theorem 2.2(2), under the present hypotheses. This result will prove to be of crucial importance in our approach of a particular case of Chinburg's $\Omega_3$–Conjecture (see §4).

Let us first remark that, under the present hypotheses,

$$r_F = \begin{cases} 1, & \text{for } F \neq k \\ 0, & \text{for } F = k \end{cases}, \quad r_\psi = \begin{cases} 1, & \text{for } \psi \neq \mathbf{1}_G \\ 0, & \text{for } \psi = \mathbf{1}_G, \end{cases} \tag{12}$$

for any $F$ and $\psi$ as in §1. For any intermediate field $F$, let $U_{F,S}$ be the group of $S$–units of $F$ ($U_S = U_{K,S}$), let $\mu_F$ be the group of roots of unity in $F^\times$, let $\overline{U}_{F,S} \stackrel{\text{def}}{=} U_{F,S}/\mu_K$ and $\overline{\mathcal{E}}_S \stackrel{\text{def}}{=} \mathcal{E}_S/\mu_K$, and let $a_F = \mathrm{Ann}_{\mathbf{Z}[G(F/k)]}(\mu_F)$. As a consequence of $|S| = 1$, we obviously have the equalities $U_{k,S} = \mu_k$ and $\overline{U}_{k,S} = \{1\}$. We also have the following characterization of $a_F$ (see [21], Chpt.IV, §1, Lemma 1):

**Lemma 3.1 (Tate).** *The ideal $a_F$ is generated over $\mathbf{Z}[G(F/k)]$ by the set*

$$\{1 - \sigma_v^{-1} \cdot Nv |\ v \text{ prime in } k,\ v \notin S'\},$$

*for any finite set of primes $S'$ of $k$, containing all the primes which ramify in $F/k$.*

An immediate consequence of Lemma 3.1 is that, if $F$ is any intermediate field, and if $\pi_F : \mathbf{Z}[G] \longrightarrow \mathbf{Z}[G(F/k)]$ is the natural projection given by $\pi_F(\sigma) = \sigma|_F$, for all $\sigma \in G$, then

$$\pi_F(a_K) = a_F. \tag{13}$$

We call a pair $(F, S')$ as in Definition 1.1 $S$–admissible if $F \neq k$, and for every finite set $T$ of primes in $k$, such that $S' \cap T = \emptyset$, the pair $(F, S')$ is $(S, T)$–admissible, in the sense of Definition 1.1. For any field $L$ of characteristic 0, and any $\psi \in \widehat{G}(L)$, $\psi \neq \mathbf{1}_G$, a pair $(F, S')$ is called $(S, \psi)$–admissible if, for every $T$ as above, the pair $(F, S')$ is $(S, T, \psi)$–admissible, in the sense of Definition 1.5. An obviously equivalent and certainly easier to handle definition of these concepts is the following:

**Definition 3.2.** A pair $(F, S')$ consisting of an intermediate field $F$, $F \neq k$, and a finite set $S'$ of primes in $k$, is called $S$–admissible if it satisfies the following conditions:

(1) $v_0 \in S'$
(2) $|S'| \geq 3$
(3) $S'$ contains all the primes which ramify in $F/k$.

If $\psi$ is as above and if, in addition to (1)—(3), $(F, S')$ satisfies

(4) $r_{\psi, S'} = 1$,

then $(F, S')$ is called $(S, \psi)$–admissible.

For an $S$–admissible pair $(F, S')$, let

$$\Theta_{F,S'}(s) = \sum_{\chi \in \widehat{G(F/k)}} L_{S'}(s, \chi) \cdot e_{\chi^{-1}}$$

be the $S'$–Stickelberger function associated to $F/k$, where $L_{S'}(s, \chi)$ is the Artin $L$–function associated to $\chi$, with Euler factors at primes in $S'$ removed.

Let us fix a prime $w_0 | v_0$ in K. In what follows we will denote by $w_0$ the prime lying below this fixed prime in any intermediate field $F$ as well. For an $S$–admissible pair $(F, S')$ we will denote by $R_{w_0}$ the regulator map

$$\mathbf{C}U_{F,S',T} = \mathbf{C}U_{F,S'} \xrightarrow{R_W} \mathbf{C}\left[G\left(F/k\right)\right]$$

associated as in §0 to $W = \{w_0\}$ and $r = 1$. The following is a reformulation in this context of Brumer–Stark's Conjecture proved independently by Deligne (see [21], Chpt.V) and Hayes (see [10]):

**Theorem 3.3 (Deligne, Hayes).** *For any $S$–admissible pair $(F, S')$ there exists a unique element $\eta_{F,S'} \in \left(\overline{U}_{F,S}\right)_{1, S'}$ satisfying the following properties:*

(1) $\dfrac{1}{e_F} R_{w_0}(\eta_{F,S'}) = \Theta'_{F,S'}(0)$

(2) $\left[\prod_{v \in T} \left(1 - \sigma_v^{-1} \cdot Nv\right)\right] \cdot \left(\dfrac{1}{e_F} \eta_{F,S'}\right) \in (U_{F,S,T})_{1,S'}$ *in* $\mathbf{Q}U_{F,S,T}$, *for every finite, nonempty set $T$ of primes in $k$, such that $T \cap S' = \emptyset$.*

Let $(F, S')$ and $T$ be as in the Theorem above. Rubin's Conjecture (Conjecture B in [15, §1.6]) applied to the set of data $(F, S', T, 1)$ predicts the existence of a unique element $\varepsilon_{F,S',T} \in (U_{F,S',T})_{1,S'}$ (see Lemma 1.5.3(2) in [15] as well) satisfying the regulator condition

$$R_{w_0}(\varepsilon_{F,S',T}) = \Theta'_{F,S',T}(0) \,.$$

Due to the obvious relation

$$\Theta'_{F,S',T}(0) = \left[\prod_{v \in T} \left(1 - \sigma_v^{-1} \cdot Nv\right)\right] \cdot \Theta'_{F,S'}(0)$$

and to the uniqueness of $\eta_{F,S'}$ and $\varepsilon_{F,S',T}$, the link between Theorem 3.3 above and Rubin's Conjecture in this context is given by

$$\varepsilon_{F,S',T} = \left[\prod_{v \in T} \left(1 - \sigma_v^{-1} \cdot Nv\right)\right] \cdot \left(\dfrac{1}{e_F} \eta_{F,S'}\right) \,. \tag{14}$$

This relation shows that Theorem 3.3 implies Rubin's Conjecture for $r = 1$. One can actually show that they are equivalent (see [17]).

**Lemma 3.4.** *Let $\ell$ be a prime number, such that $\gcd(\ell, g) = 1$, let $\psi \in \widehat{G}(\mathbf{Q}_\ell)$, $\psi \neq \mathbf{1}_G$, and let $\alpha_\psi \in a_K \otimes \mathbf{Z}_\ell$ such that $(a_K \otimes \mathbf{Z}_\ell)^\psi = \alpha_\psi D_\psi$. (Such an element always exists because $D_\psi$ is a discrete valuation ring.) Then $(\overline{\mathcal{E}}_S \otimes \mathbf{Z}_\ell)^\psi$ is generated as a $D_\psi$-module by*

$$\alpha_\psi \cdot \left[\frac{1}{e_{K_\psi}} \eta_{K_\psi, S'}\right],$$

*for any $S'$ such that the pair $(K_\psi, S')$ is $(S, \psi)$–admissible.*

**Proof.** Let us fix a set $S'$ such that the pair $(K_\psi, S')$ is $(S, \psi)$–admissible. According to Lemma 1.8(1) (for $r_\psi = 1$) and to the definition of $\mathcal{E}_S$, we have the following equality

$$\left(\overline{\mathcal{E}}_S \otimes \mathbf{Z}_\ell\right)^\psi = \sum_T D_\psi \cdot \varepsilon_{K_\psi, S', T}, \tag{15}$$

where the sum is taken with respect to all finite, nonempty sets $T$ of primes in $k$, such that $T \cap S' = \emptyset$. Lemma 3.1 shows that for any such $T$, $\prod_{v \in T} (1 - \sigma_v^{-1} \cdot Nv) \in a_{K_\psi}$, and therefore (14) and (15) show that

$$\left(\overline{\mathcal{E}}_S \otimes \mathbf{Z}_\ell\right)^\psi \subseteq D_\psi \left[\alpha_\psi \cdot \frac{1}{e_{K_\psi}} \eta_{K_\psi, S'}\right].$$

On the other hand, Lemma 3.1 also shows that one can write

$$\alpha_\psi = \sum_{v \notin S'} a_v \cdot \left(1 - \sigma_v^{-1} Nv\right),$$

with $a_v \in D_\psi$, almost all equal to 0. Let $T_v = \{v\}$, for any $v \notin S'$. The equality above, combined with (14) and the definition of $\mathcal{E}_S$, gives

$$\alpha_\psi \cdot \frac{1}{e_{K_\psi}} \eta_{K_\psi, S'} = \sum_{v \notin S'} a_v \cdot \left(1 - \sigma_v^{-1} Nv\right) \left(\frac{1}{e_{K_\psi}} \eta_{K_\psi, S'}\right) =$$

$$= \sum_{v \notin S'} a_v \cdot \varepsilon_{K_\psi, S', T_v} \in \left(\overline{\mathcal{E}}_S \otimes \mathbf{Z}_\ell\right)^\psi.$$

This concludes the proof of Lemma 3.4. □

**Lemma 3.5.** *For any $\ell$, $\psi$ and $S'$ as in Lemma 3.4, the natural morphisms of $\mathbf{Z}[G]$–modules involved in the following commutative diagrams:*

$$\begin{array}{ccc} \overline{U}_S & \longrightarrow & \overline{U}_{S'} \\ \uparrow & & \uparrow \\ \overline{U}_{K_\psi, S} & \longrightarrow & \overline{U}_{K_\psi, S'} \end{array} \qquad \begin{array}{ccc} A_S & \longrightarrow & A_{S'} \\ \uparrow & & \uparrow \\ A_{K_\psi, S} & \longrightarrow & A_{K_\psi, S'} \end{array}$$

*become $D_\psi$–isomorphisms, when tensored with $D_\psi$ over $\mathbf{Z}[G]$.*

**Proof.** Same as the proof of Lemma 1.7. □

**Proposition 3.6.** *Let $\ell$, $\psi$ and $S'$ be as in Lemma 3.4. Then*

$$|(U_S/\mathcal{E}_S \otimes \mathbf{Z}_\ell)^\psi| \geq |(A_S \otimes \mathbf{Z}_\ell)^\psi|.$$

**Proof.** Let $T$ be a finite, nonempty set of primes in $k$, such that $T \cap S' = \emptyset$. After tensoring with $D_\psi$ over $\mathbf{Z}[1/g][G]$, Theorem 0 (2) applied to the set of data $(K_\psi, S', T, 1)$ gives us the inclusion

$$D_\psi \cdot \varepsilon_{K_\psi, S', T} \subseteq \mathrm{Fitt}_{D_\psi}\left(\left(A_{K_\psi, S', T} \otimes \mathbf{Z}_\ell\right)^\psi\right) \cdot \left(U_{K_\psi, S', T} \otimes \mathbf{Z}_\ell\right)^\psi. \qquad (16)$$

On the other hand, after tensoring with $D_\psi$ over $\mathbf{Z}[G]$, the surjective $\mathbf{Z}[G]$–morphism

$$A_{K_\psi, S', T} \twoheadrightarrow A_{K_\psi, S'}$$

(see [15, (1)]) gives a surjective morphism of $D_\psi$–modules:

$$\left(A_{K_\psi, S', T} \otimes \mathbf{Z}_\ell\right)^\psi \twoheadrightarrow \left(A_{K_\psi, S'} \otimes \mathbf{Z}_\ell\right)^\psi$$

which, according to [15, §1.4(e)], gives the following inclusion at the level of Fitting ideals

$$\mathrm{Fitt}_{D_\psi}\left(A_{K_\psi, S', T} \otimes \mathbf{Z}_\ell\right)^\psi \subseteq \mathrm{Fitt}_{D_\psi}\left(A_{K_\psi, S} \otimes \mathbf{Z}_\ell\right)^\psi. \qquad (17)$$

If we now combine inclusions (16) and (17) with the obvious $U_{S', T} \subseteq \overline{U}_{S'}$ and with Lemma 3.5, we obtain

$$D_\psi \cdot \varepsilon_{K_\psi, S', T} \subseteq \mathrm{Fitt}_{D_\psi}\left(A_S \otimes \mathbf{Z}_\ell\right)^\psi \cdot \left(\overline{U}_S \otimes \mathbf{Z}_\ell\right)^\psi,$$

for every $T$ as above. This last relation together with (15) show that

$$\left(\overline{\mathcal{E}}_S \otimes \mathbf{Z}_\ell\right)^\psi \subseteq \mathrm{Fitt}_{D_\psi}\left(A_S \otimes \mathbf{Z}_\ell\right)^\psi \cdot \left(\overline{U}_S \otimes \mathbf{Z}_\ell\right)^\psi. \qquad (18)$$

On the other hand, since $r_{\psi, S} = 1$, we have an isomorphism of $D_\psi$–modules $\left(\overline{U}_S \otimes \mathbf{Z}_\ell\right)^\psi \xrightarrow{\sim} D_\psi$. Relation (18), together with [15, §1.4(f)], therefore imply that

$$|(U_S/\mathcal{E}_S \otimes \mathbf{Z}_\ell)^\psi| = |(\overline{U}_S/\overline{\mathcal{E}}_S \otimes \mathbf{Z}_\ell)^\psi| \geq$$
$$\geq \left[D_\psi : \mathrm{Fitt}_{D_\psi}\left(A_S \otimes \mathbf{Z}_\ell\right)^\psi\right] =$$
$$= |(A_S \otimes \mathbf{Z}_\ell)^\psi|,$$

which concludes the proof of Proposition 3.6. $\square$

Our next goal is to prove the following:

**Proposition 3.7.** *Let $\ell$ be a prime number such that $\gcd(\ell, g) = 1$. Then*

$$\prod_{\psi \neq \mathbf{1}_G} |(U_S/\mathcal{E}_S \otimes \mathbf{Z}_\ell)^\psi| = \prod_{\psi \neq \mathbf{1}_G} |(A_S \otimes \mathbf{Z}_\ell)^\psi|,$$

*where the products above are taken with respect to characters $\psi \in \widehat{G}(\mathbf{Q}_\ell)$.*

In order to prove the statement above we employ techniques similar to those used by Rubin in [16]. Let us fix a prime number $\ell$, $\gcd(\ell, g) = 1$, a completion $\mathbf{C}_\ell$ of the algebraic closure of $\mathbf{Q}_\ell$, with respect to the normalized $\ell$–adic valuation $v$, an embedding $\mathbf{C}_\ell \hookrightarrow \mathbf{C}$, and a finite, unramified extension $O$ of $\mathbf{Z}_\ell$, inside $\mathbf{C}_\ell$, containing the values of all characters $\chi \in \widehat{G}(\mathbf{C}_\ell)$. In what follows we will denote by $\psi$ a generic element of $\widehat{G}(\mathbf{Q}_\ell)$, and by $\chi$ a generic element of $\widehat{G}(\mathbf{C}_\ell)$. The extension of the normalized $\ell$–adic valuation $v$ to $\mathbf{C}_\ell$ (which will also be denoted by $v$) gives an equivalence relation $\sim$ on $\mathbf{C}_\ell^*$, defined by:

$$a \sim b \text{ if } v(ab^{-1}) = 0, \quad \forall a, b \in \mathbf{C}_\ell^*.$$

Let $X_S$ be the $\mathbf{Z}[G]$–module defined by

$$X_S = \{\sum_{w|v_0} a_w \cdot w \mid a_w \in \mathbf{Z}, \sum_{w|v_0} a_w = 0\}.$$

As in [15, §1], we have an injective $\mathbf{Z}[G]$–morphism

$$\overline{U}_S \xrightarrow{\lambda_S} \mathbf{R} X_S,$$

given by $\lambda_S(\bar\alpha) = \sum_{w|v_0} -\log|\alpha|_w \cdot w$, for all $\alpha \in U_S$, which becomes a $\mathbf{C}[G]$–isomorphism

$$\mathbf{C}\overline{U}_S \xrightarrow{\lambda_S}_\sim \mathbf{C} X_S,$$

when extended by $\mathbf{C}$–linearity to $\mathbf{C}\overline{U}_S$.

Let us notice that since $S = \{v_0\}$, the definitions of $R_{w_0}$ and $\lambda_S$ show that

$$\lambda_S(u) = R_{w_0}(u) \cdot w_0, \quad \forall u \in \overline{U}_S. \tag{19}$$

Let $\mathbf{Z}_\ell[G] = \bigoplus_{\psi \in \widehat{G}(\mathbf{Q}_\ell)} D_\psi$ be the usual decomposition of $\mathbf{Z}_\ell[G]$ in a direct sum of Dedekind domains. An argument similar to the one in [15], Lemma 1.5.1, combined with the second equality in (12) shows that we have isomorphisms of $\mathbf{Z}_\ell[G]$–modules

$$\mathbf{Z}_\ell \overline{U}_S \xrightarrow{\sim} \mathbf{Z}_\ell X_S \xrightarrow{\sim} \bigoplus_{\substack{\psi \in \widehat{G}(\mathbf{Q}_\ell) \\ \psi \neq \mathbf{1}_G}} D_\psi, \tag{20}$$

and equalities

$$(OX_S)^\chi = O \cdot (e_\chi w_0), \quad \forall \chi \in \widehat{G}(\mathbf{C}_\ell), \quad \chi \neq \mathbf{1}_G. \tag{21}$$

These show in particular that one can find an injective $\mathbf{Z}[G]$–morphism

$$X_S \xrightarrow{f} \overline{U}_S,$$

which becomes a $\mathbf{Z}_\ell[G]$–isomorphism, when tensored with $\mathbf{Z}_\ell$.

Let $\widehat{R}_S$, $R_{S,\chi}$ and $R_{K,S}$ be the determinants of $(\lambda_S \circ f)$, $(\lambda_S \circ f)^\chi$ and $\lambda_S$, respectively, computed with respect to $O$–bases of $OX_S$, $(OX_S)^\chi$, and $\mathbf{Z}$–bases of $X_S$ and $\overline{U}_S$ respectively, for every $\chi \in \widehat{G}(\mathbf{C}_\ell)$. For every $\psi \in \widehat{G}(\mathbf{Q}_\ell)$ we set $R_{S,\psi} \stackrel{\text{def}}{=} \prod_{\chi|\psi} R_{S,\chi}$. As a consequence of $\gcd(\#\mathrm{coker}(f), \ell) = 1$, we have the following $\ell$–adic equivalences:

$$R_{K,S} \sim \widehat{R}_S \sim \prod_{\chi \in \widehat{G}(\mathbf{C}_\ell)} R_{S,\chi} = \prod_{\psi \in \widehat{G}(\mathbf{Q}_\ell)} R_{S,\psi}. \tag{22}$$

**Lemma 3.8.** *Let $\psi \in \widehat{G}(\mathbf{Q}_\ell)$, $\psi \neq \mathbf{1}_G$, and $u \in \mathbf{Z}_\ell \overline{U}_S$, such that the index $\left[(\mathbf{Z}_\ell \overline{U}_S)^\psi : D_\psi u\right]$ is finite. Then*

$$\left[(\mathbf{Z}_\ell \overline{U}_S)^\psi : D_\psi u\right] \sim R_{S,\psi}^{-1} \cdot \left(\prod_{\chi|\psi} R_{w_0,\chi}(u)\right),$$

*where $R_{w_0,\chi} \stackrel{\text{def}}{=} \chi \circ R_{w_0}$.*

**Proof.** As a direct consequence of (19), (20) and (21), one can easily prove that

$$\left[(O\overline{U}_S)^\chi : O(e_\chi u)\right] = \ell^{v\left(R_{S,\chi}^{-1} \cdot R_{w_0,\chi}(u)\right) \cdot [O:\mathbf{Z}_\ell]}, \quad \forall \chi|\psi.$$

By taking the product of the equalities above, for all $\chi|\psi$, one obtains the statement in the lemma. $\square$

**Corollary 3.9.** *Let $\ell$, $\psi$ and $\alpha_\psi$ be as in Lemma 5.3.4. Then*

$$|(U_S/\mathcal{E}_S \otimes \mathbf{Z}_\ell)^\psi| = R_{S,\psi}^{-1} \cdot \prod_{\chi|\psi} R_{w_0,\chi}\left(\alpha_\psi \cdot \frac{1}{e_{K_\psi}} \eta_{K_\psi, S_\psi}\right).$$

**Proof.** Immediate consequence of Theorem 2.2(1), Lemma 3.4 and Lemma 3.8. $\square$

Let $\chi \in \widehat{G}(\mathbf{C})$. (Under the chosen embeding $\mathbf{C}_\ell \hookrightarrow \mathbf{C}$, $\chi$ can be viewed as an element of $\widehat{G}(\mathbf{C}_\ell)$ as well.) Let $L_S(s, \chi)$ be the Artin $L$–function associated to $\chi$, with the Euler factor at $v_0$ removed:

$$L_S(s, \chi) = \left(1 - \chi(\sigma_{v_0}) Nv_0^{-s}\right) \cdot L(s, \chi),$$

where $L(s, \chi)$ is the "complete" Artin $L$–function associated to $\chi$. For every $\psi \in \widehat{G}(\mathbf{Q}_\ell)$, $\psi \neq \mathbf{1}_G$, let $S_\psi$ be a fixed finite set of primes in $k$, such that the pair $(K_\psi, S_\psi)$ is $(S, \psi)$–admissible. Then every $\chi \in \widehat{G}(\mathbf{C}_\ell)$, $\chi | \psi$, gives a relation

$$L_{S_\psi}(s, \chi) = \prod_{v \in S_\psi \setminus S'_\psi} \left(1 - \chi(\sigma_v) Nv^{-s}\right) \cdot L_S(s, \chi),$$

where $S'_\psi = \{v \in S_\psi |\, v \text{ ramifies in } K_\psi/k, \text{ or } v = v_0\}$. Due to the fact that $(K_\psi, S_\psi)$ is $(S, \psi)$–admissible, we have $r_{\psi, S_\psi} = 1$, and therefore $\chi(\sigma_v)$ is a root of unity of order $g$, different from 1, for all $v \in S_\psi \setminus S'_\psi$ and all $\chi | \psi$. Since $\gcd(\ell, g) = 1$, this shows that for $\chi$ and $v$ as above we have $(1 - \chi(\sigma_v)) \in O^\times$. We therefore have the following $\ell$–adic equivalences

$$L'_{S_\psi}(0, \chi) = \prod_{v \in S_\psi \setminus S'_\psi} (1 - \chi(\sigma_v)) \cdot L'_S(0, \chi) \sim L'_S(0, \chi), \tag{23}$$

for every $\psi$ and $\chi$ as above. We are now prepared to prove Proposition 3.7.

**Proof of Proposition 3.7.** Let $\zeta_{K,S}(s)$ and $\zeta_{k,S}(s)$ be the $S$–zeta functions of $K$ and $k$ respectively. Let $h_{K,S} = |A_S|$, $h_{k,S} = |A_{k,S}|$, and for every $\psi \in \widehat{G}(\mathbf{Q}_\ell)$, let $h^\psi_{K,S} = |(A_S \otimes \mathbf{Z}_\ell)^\psi|$ and $e^\psi_K = |(\mu_K \otimes \mathbf{Z}_\ell)^\psi|$. We have the equality

$$\frac{\zeta_{K,S}(s)}{\zeta_{k,S}(s)} = \prod_{\psi \neq \mathbf{1}_G} \left( \prod_{\chi | \psi} L_S(\chi, s) \right). \tag{24}$$

Since $S$ consists of a single prime which splits completely in $K/k$, we have the following equalities (see [21], Chpt.I):

$$\operatorname{ord}_{s=0} \zeta_{k,S}(s) = 0\,,\ \operatorname{ord}_{s=0} \zeta_{K,S}(s) = [K:k] - 1\,,\ \operatorname{ord}_{s=0} L_S(\chi, s) = 1\,, \forall \chi \neq \mathbf{1}_G.$$

The $S$–class number formula (see [21], Chpt.I) together with equality (24) therefore give us

$$\frac{h_{K,S}}{h_{k,S}} \cdot \frac{R_{K,S}}{R_{k,S}} \cdot \frac{e_K}{e_k} = \prod_{\psi \neq \mathbf{1}_G} \left( \prod_{\chi | \psi} L'_S(\chi, 0) \right). \tag{25}$$

According to (22), (25), and Theorem 3.3(1) applied to the pairs $(K_\psi, S_\psi)$, for all $\psi \in \widehat{G}(\mathbf{Q}_\ell)$, $\psi \neq \mathbf{1}_G$, we therefore have

$$\prod_{\psi \neq \mathbf{1}_G} h_{K,S}^\psi \cdot \prod_{\psi \neq \mathbf{1}_G} R_{S,\psi} \cdot \left( \prod_{\psi \neq \mathbf{1}_G} e_K^\psi \right)^{-1} \sim \prod_{\psi \neq \mathbf{1}_G} \left( \prod_{\chi | \psi} R_{w_0, \chi} \left( \frac{1}{e_{K_\psi}} \eta_{K_\psi, S_\psi} \right) \right). \quad (26)$$

Let us choose $\alpha_\psi \in (a_K \otimes \mathbf{Z}_\ell)$ as in Lemma 3.4, for every $\psi$ as above. The flatness of $D_\psi$ and the cyclicity of $\mu_K$ as $\mathbf{Z}[G]$–modules give us the following exact sequence of $D_\psi$–modules

$$0 \longrightarrow (a_K \otimes \mathbf{Z}_\ell)^\psi \longrightarrow D_\psi \longrightarrow (\mu_K \otimes \mathbf{Z}_\ell)^\psi \longrightarrow 0,$$

for every $\psi$. This together with (13) give

$$e_{K_\psi} = e_K^\psi = [D_\psi : \alpha_\psi D_\psi]. \quad (27)$$

Let $\psi$ as above, and for a fixed $\chi_0 | \psi$, let $\mathbf{Z}_\ell[\chi_0]$ be the discrete valuation domain obtained by adjoining the values of $\chi_0$ to $\mathbf{Z}_\ell$. Then we have the equalities

$$\begin{aligned}
[D_\psi : \alpha_\psi D_\psi] &= [\mathbf{Z}_\ell[\chi_0] : (\chi_0(\alpha_\psi))] = \\
&= \ell^{v(\chi_0(\alpha_\psi)) \cdot [\mathbf{Z}_\ell[\chi_0] : \mathbf{Z}_\ell]} = \\
&= \ell^{v\left( \prod_{\chi | \psi} \chi(\alpha_\psi) \right)}.
\end{aligned} \quad (28)$$

Equalities (27) and (28) show that we have an $\ell$–adic equivalence

$$e_{K_\psi} \sim \prod_{\chi | \psi} \chi(\alpha_\psi),$$

for every $\psi \neq \mathbf{1}_G$. This equivalence together with the $G$–linearity of $R_{w_0}$ show that, if we multiply (26) by $\prod_{\psi \neq \mathbf{1}_G} R_{S,\psi}^{-1} \left( \prod_{\chi | \psi} \chi(\alpha_\psi) \right)$, we obtain

$$\prod_{\psi \neq \mathbf{1}_G} h_{K,S}^\psi \sim \prod_{\psi \neq \mathbf{1}_G} \left[ R_{S,\psi}^{-1} \cdot \prod_{\chi | \psi} R_{w_0, \chi} \left( \alpha_\psi \cdot \frac{1}{e_{K_\psi}} \eta_{K_\psi, S_\psi} \right) \right]. \quad (29)$$

If we now combine (29) with Corollary 3.9, we obtain

$$\prod_{\psi \neq \mathbf{1}_G} h_{K,S}^\psi \sim \prod_{\psi \neq \mathbf{1}_G} |(U_S / \mathcal{E}_S \otimes \mathbf{Z}_\ell)^\psi|,$$

which concludes the proof of Proposition 3.7. $\square$

By combining Propositions 3.6 and 3.7 we obtain the desired result

**Theorem 3.10.** *Let $S$ be a set consisting of a single prime of $k$, which splits completely in $K/k$. Then:*

(1) *The index $[U_S : \mathcal{E}_S]$ is finite.*
(2) *For every prime number $\ell$ such that $\gcd(\ell, g) = 1$, and every $\psi \in \widehat{G}(\mathbf{Q}_\ell)$, $\psi \neq \mathbf{1}_G$,*
$$|(U_S/\mathcal{E}_S \otimes \mathbf{Z}_\ell)^\psi| = |(A_S \otimes \mathbf{Z}_\ell)^\psi|.$$

**Remark.** Equalities similar to the ones in Theorems 1.4(2), 2.2(2) and 3.10(2) were proved in the number field case, for $k = \mathbf{Q}$, as consequences of the Main Conjecture in Iwasawa Theory (see [12]), or by using Kolyvagin–Euler System techniques (see[11], Appendix). It is conceivable, as Rubin shows in [17, §6], that the elements $\Phi(\varepsilon_{F,S',T})$ give Kolyvagin–Euler Systems in the function field case as well. One could therefore hope to be able to prove the theorems above by using Kolyvagin–Euler System techniques. However, there would be a major drawback: these techniques would fail to give information at the prime $p = \mathrm{char}(k)$. The methods we are using do not force us to stay away from the prime $p$. This is mainly due to the fact that we are not avoiding $p$ in Theorem 0. The cause for this lies deeper, in the links between crystalline and $p$–adic étale cohomology theories, and in the fact that the action of the geometric Frobenius morphism on the crystalline cohomology groups gives the right $L$–functions (see Appendix of [15] or [14]). The idea of using crystalline cohomology in order to understand the structure of the $p$–part of class groups of function fields has been used in the past by Goss and Sinnott (see [7]).

## 4. On Chinburg's $\Omega_3$–Conjecture for function fields

4.1. THE CONJECTURE

Let $K/k$ be any finite, Galois extension of global fields (number fields or function fields of characteristic $p > 0$.) Let $G = G(K/k)$, and let $\mathrm{K}_0(\mathbf{Z}[G])$ be the Grothendieck group of isomorphism classes of finitely generated, projective $\mathbf{Z}[G]$–modules. If $P$ is a finitely generated, projective $\mathbf{Z}[G]$–module, then $P \otimes_\mathbf{Z} \mathbf{Q}$ is $\mathbf{Q}[G]$–free (see [18]) and we define
$$\mathrm{rank}(P) \stackrel{\mathrm{def}}{=} \mathrm{rank}_{\mathbf{Q}[G]}(P \otimes_\mathbf{Z} \mathbf{Q}).$$

We obtain this way a surjective group morphism

$$\mathrm{K}_0(\mathbf{Z}[G]) \xrightarrow{\mathrm{rank}} \mathbf{Z},$$

defined by $\mathrm{rank}([P]) = \mathrm{rank}(P)$, for the class $[P]$ in $\mathrm{K}_0(\mathbf{Z}[G])$ of any finitely generated, projective $\mathbf{Z}[G]$–module $P$.

**Definition 4.1.1.** The projective class–group $\mathrm{Cl}(\mathbf{Z}[G])$ of $\mathbf{Z}[G]$ is defined by

$$\mathrm{Cl}(\mathbf{Z}[G]) = \ker\left(\mathrm{K}_0(\mathbf{Z}[G]) \xrightarrow{\mathrm{rank}} \mathbf{Z}\right).$$

In [6] Fröhlich defined an invariant of *analytic* nature $\mathrm{W}_{K/k} \in \mathrm{Cl}(\mathbf{Z}[G])$ for finite, Galois extensions $K/k$ of number fields, by means of the Artin root numbers associated to the irreducible, symplectic representations of $G$. In [5] Chinburg shows that the definition of $\mathrm{W}_{K/k}$ can be carried through in the function field case as well. In particular, if the group $G$ has no irreducible, symplectic representations (as it is the case when $G$ is abelian, for example), then $\mathrm{W}_{K/k} = 0$ in $\mathrm{Cl}(\mathbf{Z}[G])$.

In [3] and [4] Chinburg defined an invariant of *arithmetic* nature $\Omega(K/k, 3) \in \mathrm{Cl}(\mathbf{Z}[G])$ for finite, Galois extensions of number fields, measuring the complexity of the $\mathbf{Z}[G]$–module structure of the group $U_S$ of $S$–units in $K$, for a sufficiently large finite set $S$ of primes in $k$. In [5] he extends this definition to the function field case.

The link between the invariants described above is believed to be the following (see [3], [4] and [5]):

**Conjecture (Chinburg).** *For any finite, Galois extension of global fields*

$$\Omega(K/k, 3) = W_{K/k}.$$

Chinburg shows in [3, §IX] how, in the case of a cyclic extension of prime degree of the field $\mathbf{Q}$ of rational numbers, the conjecture above follows from Gras' Conjecture satisfied by the group of cyclotomic units in that context. Following Chinburg's ideas, we are going to show in the next section that the conjecture above holds true for cyclic extensions of prime degree of function fields, as a consequence of Theorem 3.10 (our analogue of Gras' Conjecture in the function field setting). Chinburg's Conjecture in this situation follows from work of S. Bae [1] as well. However Bae's approach is completely different from ours.

## 4.2. CYCLIC EXTENSIONS OF PRIME DEGREE

Throughout this section we are going to assume $K/k$ to be a cyclic function field extension of prime degree. Since in this case there are no irreducible, symplectic representations of $G$, we have

$$W_{K/k} = 0, \tag{29}$$

and therefore Chinburg's Conjecture asserts that $\Omega(3, K/k) = 0$ in $\mathrm{Cl}(\mathbf{Z}[G])$. Since $|G| = 2$ implies $\mathrm{Cl}(\mathbf{Z}[G]) = 0$ (see [3, §IX]), we can assume from now on that $|G| > 2$.

Let $G_0(\mathbf{Z}[G])$ be the Grothendieck group of isomorphism classes of finitely generated $\mathbf{Z}[G]$–modules, and let

$$h : K_0(\mathbf{Z}[G]) \longrightarrow G_0(\mathbf{Z}[G])$$

be the Cartan morphism, taking the class $[P]$ of a finitely generated, projective $\mathbf{Z}[G]$–module $P$ in $K_0(\mathbf{Z}[G])$, into its class $(P)$ in $G_0(\mathbf{Z}[G])$.

For a finite set $S$ of primes of $k$, let $S_K$ be the set of primes of $K$ lying above primes in $S$. Let us consider the following $\mathbf{Z}[G]$–modules:

$$Y_S = \{\sum_{w \in S_K} a_w \cdot w \mid a_w \in \mathbf{Z}\}, \quad X_S = \{\sum_{w \in S_K} a_w \cdot w \mid a_w \in \mathbf{Z}, \sum_{w \in S_K} a_w = 0\}.$$

The following statements are proved in [5] and [19] respectively:

**Proposition 4.2.1 (Chinburg).** *The image of $\Omega(3, K/k)$ via the Cartan map $h$ satisfies*

$$h(\Omega(3, K/k)) = (U_S) - (X_S) - (A_S),$$

*for any finite, sufficiently large set $S$ of primes in $k$.*

**Proposition 4.2.2 (Reiner).** *If $G$ is a cyclic group of prime order then, the Cartan map*

$$h : K_0(\mathbf{Z}[G]) \longrightarrow G_0(\mathbf{Z}[G])$$

*is injective.*

According to these Propositions, all we have to do in order to prove the Conjecture above under the present assumptions, is show that $(U_S) - (X_S) - (A_S) = 0$ in $G_0(\mathbf{Z}[G])$, for a sufficiently large set of primes $S$. Although one can make the expression "sufficiently large" very precise (see [5] ), the next result shows that there is no need for that in the present situation.

**Lemma 4.2.3.** *The class*

$$c_S \stackrel{\mathrm{def}}{=} (U_S) - (X_S) - (A_S) \in G_0(\mathbf{Z}[G])$$

*does not depend on the finite, nonempty set $S$ of primes in $k$.*

**Proof.** There is an exact sequence of $\mathbf{Z}[G]$–modules

$$0 \longrightarrow X_S \longrightarrow Y_S \stackrel{s}{\longrightarrow} \mathbf{Z} \longrightarrow 0,$$

where $s\left(\sum_{w \in S_K} a_w \cdot w\right) \stackrel{\mathrm{def}}{=} \sum_{w \in S_K} a_w$, for all $\sum_{w \in S_K} a_w \cdot w \in Y_S$. This gives the following relation in $G_0(\mathbf{Z}[G])$

$$(X_S) = (Y_S) - (\mathbf{Z}) . \tag{30}$$

In order to prove the desired statement, it is obviously enough to prove that, if $S_1$ and $S_2$ are two sets of primes as in the statement of the lemma, satisfying $S_1 \subseteq S_2$, then $c_{S_1} = c_{S_2}$ in $G_0(\mathbf{Z}[G])$.

Let us fix such sets $S_1$ and $S_2$. Let $A_{S_1,S_2}$ be the $\mathbf{Z}[G]$–submodule of $A_{S_1}$ generated by the classes of the ideals in $S_{1,2} \stackrel{\mathrm{def}}{=} S_{2,K} \setminus S_{1,K}$, and let $Y_{S_1,S_2}$ be the free abelian group generated by the ideals in $S_{1,2}$. We have the following exact sequences of finitely generated $\mathbf{Z}[G]$–modules:

$$0 \longrightarrow Y_{S_1} \longrightarrow Y_{S_2} \stackrel{\alpha}{\longrightarrow} Y_{S_1,S_2} \longrightarrow 0$$

$$0 \longrightarrow U_{S_1} \longrightarrow U_{S_2} \stackrel{\beta}{\longrightarrow} Y_{S_1,S_2} \stackrel{\gamma}{\longrightarrow} A_{S_1,S_2} \longrightarrow 0$$

$$0 \longrightarrow A_{S_1,S_2} \longrightarrow A_{S_1} \longrightarrow A_{S_2} \longrightarrow 0,$$

where $\alpha\left(\sum_{w \in S_{2,K}} a_w \cdot w\right) = \sum_{w \in S_{1,2}} a_w \cdot w$, $\beta(u) = \sum_{w \in S_{1,2}} \mathrm{ord}_w(u) \cdot w$, for all $u \in U_{S_2}$, and $\gamma(w) = \widehat{w}$ (the class of $w$ in $A_{S_1,S_2}$), for all $w \in S_{1,2}$. If we write down the relations given in $G_0(\mathbf{Z}[G])$ by the three exact sequences above, and eliminate the classes $(Y_{S_1,S_2})$ and $(A_{S_1,S_2})$, we obtain

$$(Y_{S_2}) - (U_{S_2}) - (A_{S_2}) = (Y_{S_1}) - (U_{S_1}) - (A_{S_1}).$$

This relation combined with (30) for $S = S_1$ and $S = S_2$ respectively, gives the desired $c_{S_1} = c_{S_2}$ in $G_0(\mathbf{Z}[G])$. $\square$

Lemma 4.2.3 shows that, in order to prove Chinburg's Conjecture under the present assumptions, we have to find a finite, nonempty set $S$ of primes in $k$, satisfying $c_S = 0$.

Let $S = \{v_0, \ldots, v_s\}$ be a fixed set of primes in $k$ satisfying the following properties:

(1) $v_0$ splits completely in $K/k$.
(2) $|S| \geq 3$
(3) $S$ contains all the primes which ramify in $K/k$.
(4) $G_{v_i} = G$, for all $i = 1, \ldots, s$.

We intend to prove that $S$ satisfies $c_S = 0$. In order to do this, we return to some of the notations and techniques employed in §3.

Let $S_0 = \{v_0\}$. Properties (1)—(3) satisfied by $S$ show that $(K, S)$ is an $S_0$–admissible pair. Let $\eta_{K,S} \in (\overline{U}_{K,S_0})_{1,S}$ be the unique element provided by Theorem 3.3 (with $F = K$, $S = S_0$ and $S' = S$). Let $\overline{E}_S$ be the $\mathbf{Z}[G]$–submodule of $\overline{U}_S = U_S/\mu_K$ defined by

$$\overline{E}_S = \left(a_K \cdot \frac{1}{e_K} \eta_{K,S}\right) \bigoplus (h_{k,S} \mathbf{Z} u_1 \oplus \mathbf{Z} u_2 \oplus \cdots \oplus \mathbf{Z} u_s),$$

where $\{u_1, \ldots, u_s\}$ is a $\mathbf{Z}$–basis of $\overline{U}_{k,S}$. (Since $\eta_{K,S} \in \left(\overline{U}_{K,S_0}\right)_{1,S}$, and $r_{\mathbf{1}_G, S} \geq 2$ (see property (2) satisfied by $S$), we have $e_{\mathbf{1}_G} \cdot \eta_{K,S} = 0$ and therefore the relation $a_K \cdot \frac{1}{e_K}\eta_{K,S} \cap \overline{U}_{k,S} = \{0\}$.)

Let $E_S \subseteq U_S$ be the preimage of $\overline{E}_S \subseteq \overline{U}_S$ under the projection

$$U_S \twoheadrightarrow \overline{U}_S = U_S/\mu_K \ .$$

Then we obviously have an exact sequence of finitely generated $\mathbf{Z}[G]$–modules

$$0 \longrightarrow \mu_K \longrightarrow E_S \longrightarrow \overline{E}_S \longrightarrow 0 \ ,$$

which gives the following relation in $\mathrm{G}_0\left(\mathbf{Z}[G]\right)$

$$(E_S) = \left(\overline{E}_S\right) + (\mu_K) \ . \tag{31}$$

**Proposition 4.2.4.**
   (1) *The index $[U_S : E_S]$ is finite.*
   (2) *For every prime number $\ell$ such that $\gcd(\ell, g) = 1$, and every $\psi \in \widehat{G}(\mathbf{Q}_\ell)$*

$$|\left(U_S/E_S \otimes \mathbf{Z}_\ell\right)^\psi| = |\left(A_S \otimes \mathbf{Z}_\ell\right)^\psi| \ .$$

**Proof.** As remarked many times before, it is enough to prove (2) above.

Let $\ell$ be as above and let us fix $\psi \in \widehat{G}(\mathbf{Q}_\ell)$, $\psi \neq \mathbf{1}_G$ for the beginning. The fact that $|G|$ is prime implies that $K_\psi = K$. Properties (1)–(4) satisfied by $S$ imply that the pair $(K = K_\psi, S)$ is $(S_0, \psi)$–admissible. Lemma 3.5 therefore shows that

$$\left(\mathbf{Z}_\ell \overline{U}_S\right)^\psi = \left(\mathbf{Z}_\ell \overline{U}_{S_0}\right)^\psi \text{ and } \left(\mathbf{Z}_\ell A_S\right)^\psi = \left(\mathbf{Z}_\ell A_{S_0}\right)^\psi \ . \tag{32}$$

The definition of $\overline{E}_S$ and Lemma 5.3.4 show that

$$\left(\mathbf{Z}_\ell \overline{E}_S\right)^\psi = \left(\mathbf{Z}_\ell a_K\right)^\psi \left(\frac{1}{e_K}\eta_{K,S}\right) = \left(\mathbf{Z}_\ell \mathcal{E}_{S_0}\right)^\psi \ . \tag{33}$$

If we now combine equalities (32) and (33) with Theorem 3.10 (for $S = S_0$) we obtain

$$|\left(U_S/E_S \otimes \mathbf{Z}_\ell\right)^\psi| = |\left(\overline{U}_S/\overline{E}_S \otimes \mathbf{Z}_\ell\right)^\psi| =$$
$$= |\left(A_S \otimes \mathbf{Z}_\ell\right)^\psi| \ ,$$

which concludes the proof of statement (2), for $\psi \neq \mathbf{1}_G$.

Let $\psi = \mathbf{1}_G$. We obviously have the following equalities:
$$\left(\mathbf{Z}_\ell \overline{E}_S\right)^{\mathbf{1}_G} = h_{k,S}\mathbf{Z}_\ell u_1 \oplus \mathbf{Z}_\ell u_2 \oplus \cdots \oplus \mathbf{Z}_\ell u_s, \quad \left(\mathbf{Z}_\ell \overline{U}_S\right)^{\mathbf{1}_G} = \mathbf{Z}_\ell \overline{U}_{k,S},$$
$$\left(\mathbf{Z}_\ell A_S\right)^{\mathbf{1}_G} = \mathbf{Z}_\ell A_{k,S}.$$

These imply that
$$|(U_S/E_S \otimes \mathbf{Z}_\ell)^{\mathbf{1}_G}| = |\left(\overline{U}_S/\overline{E}_S \otimes \mathbf{Z}_\ell\right)^{\mathbf{1}_G}| =$$
$$= |(\mathbf{Z}_\ell \otimes A_S)^{\mathbf{1}_G}|,$$

which concludes the proof of Proposition 4.2.4 (2), for $\psi = \mathbf{1}_G$ as well. $\square$

For a finite $\mathbf{Z}[G]$–module $M$, we have a $\mathbf{Z}[G]$–module direct sum decomposion
$$M = \bigoplus_\ell (M \otimes \mathbf{Z}_\ell),$$
with respect to all prime numbers $\ell$. Let $\mathcal{S}(M \otimes \mathbf{Z}_\ell)$ denote the $\mathbf{Z}[G]$– (or equivalently, the $\mathbf{Z}_\ell[G]$–) semisimplification of $M \otimes \mathbf{Z}_\ell$, for any $\ell$. We obviously have the following relation in $\mathrm{G}_0(\mathbf{Z}[G])$
$$(M) = \sum_\ell \left(\mathcal{S}(M \otimes \mathbf{Z}_\ell)\right). \tag{34}$$

On the other hand, if $\ell = |G|$ (recall that $|G|$ is a prime number), and $N$ is a simple $\mathbf{Z}[G]$–module of order a power of $\ell$, then the maximal $G$–fixed submodule $N^G$ of $N$ is nontrivial, and therefore $N^G = N$. This fact, combined again with the simplicity of $N$, shows that $|N| = \ell$. We therefore have an exact sequence of $\mathbf{Z}[G]$–modules
$$0 \longrightarrow \mathbf{Z} \xrightarrow{\times \ell} \mathbf{Z} \longrightarrow N \longrightarrow 0,$$
which shows that $(N) = 0$ in $\mathrm{G}_0(\mathbf{Z}[G])$.

This remark and (34) show that, for any finite $\mathbf{Z}[G]$–module $M$, we have the following relation in $\mathrm{G}_0(\mathbf{Z}[G])$
$$(M) = \sum_{\ell \neq |G|} \left(\mathcal{S}(M \otimes \mathbf{Z}_\ell)\right). \tag{35}$$

**Lemma 4.2.5.** *The equality*
$$(U_S) - (E_S) - (A_S) = 0$$
*holds true in* $\mathrm{G}_0(\mathbf{Z}[G])$.

**Proof.** Greenberg shows in [8, §5] that, if $G$ is a finite abelian group, $\ell$ is a prime number such that $\gcd(\ell, |G|) = 1$, and $M_1$, $M_2$ are two finite $\mathbf{Z}[G]$–modules, the following statements are equivalent:
   (1) $\mathcal{S}(M_1 \otimes \mathbf{Z}_\ell) \xrightarrow{\sim} \mathcal{S}(M_2 \otimes \mathbf{Z}_\ell)$ as $G$–modules.
   (2) $|(M_1 \otimes \mathbf{Z}_\ell)^\psi| = |(M_2 \otimes \mathbf{Z}_\ell)^\psi|$, for all $\psi \in \widehat{G}(\mathbf{Q}_\ell)$.

The statement in Lemma 4.2.5 follows now from Proposition 4.2.4, relation (35) and Greenberg's observation, with $M_1 = U_S/E_S$ and $M_2 = A_S$. $\square$

**Proposition 4.2.6.** *The equality* $(E_S) = (X_S)$ *holds true in* $\mathrm{G}_0(\mathbf{Z}[G])$.

**Proof.** Let us fix $w_i \in S_K$, $w_i | v_i$, for every $i = 0, \ldots, s$. Properties (1) and (4) satisfied by $S$ imply that we have $\mathbf{Z}[G]$–isomorphisms

$$\mathbf{Z}[G] w_0 \xrightarrow{\sim} \mathbf{Z}[G], \quad \mathbf{Z}[G] w_i \xrightarrow{\sim} \mathbf{Z}, \forall i = 1, \ldots, s,$$

with $G$ acting trivially on $\mathbf{Z}$. The definition of $Y_S$ therefore shows that

$$Y_S = \bigoplus_{0 \leq i \leq s} \mathbf{Z}[G] w_i \xrightarrow{\sim} \mathbf{Z}[G] \oplus \mathbf{Z}^s$$

as $\mathbf{Z}[G]$–modules. Relation (30) therefore shows that the following holds true in $\mathrm{G}_0(\mathbf{Z}[G])$

$$(X_S) = (\mathbf{Z}[G]) + (s - 1)(\mathbf{Z}). \tag{36}$$

We are now going to compute the class $(E_S)$ in $\mathrm{G}_0(\mathbf{Z}[G])$. According to (31) and the definition of $\overline{E}_S$, we obviously have

$$(E_S) = (\mu_K) + \left(a_K \cdot \frac{1}{e_K}\eta_{K,S}\right) + s(\mathbf{Z}). \tag{37}$$

Since $\mu_K$ is a cyclic $\mathbf{Z}[G]$–module, the definition of $a_K$ implies

$$(a_K) = (\mathbf{Z}[G]) - (\mu_K). \tag{38}$$

Proposition 4.2.4 (1) implies that

$$\mathrm{rank}_{\mathbf{Z}}\overline{E}_S = \mathrm{rank}_{\mathbf{Z}}\overline{U}_S = [K : k] + (s - 1),$$

and therefore, from the definition of $\overline{E}_S$, we have

$$\mathrm{rank}_{\mathbf{Z}}\left(a_K \cdot \frac{1}{e_K}\eta_{K,S}\right) = [K : k] - 1. \tag{39}$$

Let $\mathcal{K}$ be the $\mathbf{Z}[G]$–module defined by the exact sequence

$$0 \longrightarrow \mathcal{K} \longrightarrow a_K \xrightarrow{\pi} a_K \cdot \frac{1}{e_K}\eta_{K,S} \longrightarrow 0, \tag{40}$$

where $\pi(\alpha) = \alpha \cdot \frac{1}{e_K}\eta_{K,S}$, for all $\alpha \in a_K$. Then (39), combined with the obvious $\mathrm{rank}_{\mathbf{Z}}(a_K) = [K : k]$, shows that $\mathrm{rank}_{\mathbf{Z}}(\mathcal{K}) = 1$. But since $|G| > 2$, the only ideal of $\mathbf{Z}[G]$ of $\mathbf{Z}$–rank equal to 1 is isomorphic to $\mathbf{Z}$ with trivial $G$–action. This shows that $(\mathcal{K}) = (\mathbf{Z})$ and therefore the following equality holds true in $\mathrm{G}_0(\mathbf{Z}[G])$

$$\left(a_K \cdot \frac{1}{e_K}\eta_{K,S}\right) = (a_K) - (\mathbf{Z}).$$

Combining this equality with (37) and (38) one obtains

$$(E_S) = (\mathbf{Z}[G]) + (s - 1)(\mathbf{Z}),$$

which, according to (36), concludes the proof of Proposition 4.2.6. □

We are now prepared to prove the main result of this section:

**Theorem 4.2.7.** *If $K/k$ is a cyclic extension of prime degree of function fields of characteristic $p > 0$, then $\Omega(K/k, 3) = W_{K/k}(= 0)$ in $\mathrm{Cl}(\mathbf{Z}[G])$.*

**Proof.** Proposition 4.2.6 together with Lemma 4.2.5 show that, for the set $S$ of primes in $k$ fixed above, we have $c_S = 0$ in $\mathrm{G}_0(\mathbf{Z}[G])$. Proposition 4.2.1 and Lemma 4.2.3 thus show that $h(\Omega(K/k, 3)) = 0$ in $\mathrm{G}_0(\mathbf{Z}[G])$. The injectivity of $h$ (see Proposition 4.2.2) together with (29), therefore imply that $\Omega(K/k, 3) = 0 = W_{K/k}$ in $\mathrm{Cl}(\mathbf{Z}[G])$, which concludes the proof. $\square$

Mathematical Sciences Research Institute, 1000 Centennial Drive, Berkeley, CA, 94720–5070, USA

*E-mail address*: popescu@msri.org